\documentstyle[12pt]{article}
\newtheorem{th}{\indent Theorem.}
\newtheorem{lem}{\indent Lemma.}
\begin{document}
\begin{center}
{\large\bf On the minimal cardinality of a subset of ${\bf R}$
which is not of first category}
\par\vspace{.3cm}\par
{\bf Apoloniusz Tyszka}
\end{center}
\begin{abstract}
\def\thefootnote{}
\noindent
Let $M$ be the ideal of first category subsets of ${\bf R}$ and
non$(M) = \min \{\mbox{card } X : X \subseteq {\bf R}, X \not\in M\}$.
We consider families $\Phi$ of sequences converging to $\infty$,
with the property that for every open set $U \subseteq {\bf R}$
that is unbounded above there exists a sequence belonging to $\Phi$,
which has an infinite number of terms belonging to $U$. We present
assumptions about $\Phi$ which imply that the minimal cardinality of
$\Phi$ equals non$(M)$.
\end{abstract}
\par
We consider families $\Phi$ of sequences converging to $\infty$ that
$\Phi$ satisfies the following condition (C) defined in [5]:
\begin{itemize}
\item[(C)] {\it for every open set $U \subseteq {\bf R}$ that is
unbounded above there exists a sequence belonging to $\Phi$, which
has an infinite number of terms belonging to $U$}.
\end{itemize}
Our considerations are motivated by the following Propositions 1-2.
\par
\
{\bf Proposition 1.} Assume that: $\Phi$ is a family of sequences
converging to $\infty$, $d\in {\bf R}\cup \{-\infty,\infty\}$ and
$\Phi$ satisfies the condition (C). We claim that if
$f:{\bf R}\longrightarrow {\bf R}$ is continuous then the convergence
of $\{ f(x_n)\} _{n\in\omega}$ to $d$ for all sequences
$\{ x_n\}_{n \in \omega}$ belonging to $\Phi$ implies that
$lim_{x\rightarrow\infty} f(x)=d$.
\par
\
{\bf Proposition 2.} Assume that: $\Phi$ is a family of sequences
converging to $\infty$, $d\in {\bf R}\cup \{-\infty,\infty\}$ and
for each continuous $f:{\bf R}\longrightarrow {\bf R}$ the convergence
of $\{ f(x_n)\} _{n\in\omega}$ to $d$ for all sequences
$\{ x_n\}_{n \in \omega}$ belonging to $\Phi$ implies that
$lim_{x\rightarrow\infty} f(x)=d$. We claim that $\Phi$ satisfies the
condition (C).
\vspace{0.3cm}
\par
We present below necessary definitions and results from [2] and [3].
For the functions $f, g: \omega\longrightarrow\omega$ we define:
$$
f\leq_* g \Longleftrightarrow\{ i: f(i) > g(i)\} \ \mbox{is finite}
\ (g \ \mbox{dominates} \ f)
$$
We say that $F \subseteq\omega^\omega$ is bounded if there exists a
$g:\omega\longrightarrow\omega$ such that $f\leq_* g$ for
all $f \in F$. Let $b = \min$ \{card $F: F$ is unbounded\}.
Let $M$ be the ideal of first category subsets of ${\bf R}$ and
non$(M) = \min \{\mbox{card } X: X \subseteq {\bf R}, X \not\in M\}$.
Let $\forall^\infty$ abbreviate "for all but finitely many". We have
([1], [2] and also [4]):
$$
\begin{array}{c}
b\leq \mbox{non}(M) = \\ [.2cm]
\min \{\mbox{card } F: F\subseteq\omega^\omega
\mbox{ and } \neg \ \exists \ g\in\omega^\omega \ \forall f\in F \
\forall^\infty k \ g(k)\neq f(k)\}
\end{array}
$$
For a sequence $\{ a_n \}_{n\in\omega}$ converging to $\infty$ we
define the non-decreasing function
$f_{\{ a_n \} }:\omega\longrightarrow\omega$ (see [5]): if
$\bigcup_{n\in\omega} (a_n -1, a_n +1)\supseteq (i,\infty)$ then
$$
f_{\{ a_n \} }(i) = \max\{ j\in\omega\setminus\{ 0\}:
\bigcup_{n\in\omega}(a_n -\frac{1}{j}, a_n +\frac{1}{j}) \supseteq
(i,\infty ) \}
$$
else $f_{\{ a_n \} }(i) = 0$.
\def\theth{}
\begin{th}
If a family $\Phi$ of sequences converging to $\infty$ satisfies
condition (C) and the following conditions hold:
\begin{equation}
\{ f_{\{ a_n \} }: \{ a_n \}\in\Phi\} \mbox{ is bounded }
\end{equation}
\begin{equation}
\mbox{ if } \{ a_n\}\in\Phi \mbox{ then } \forall n\in\omega\setminus
\{ 0 \} \ 0 < a_{n+1}-a_n\leq a_n -a_{n-1}
\end{equation}
then card $\Phi\geq non(M)$.
\end{th}
\par
{\it Proof}. We start from remarks which show connections between
our theorem and results of [5]. After these remarks we present an
introductory Observation and the main Lemma.
\par
{\bf Remark 1.} ([5]) If for a family $\Phi$ of sequences converging
to $\infty$ the
family of functions $\{ f_{\{ a_n \} }: \{ a_n \}\in\Phi\}$ is
unbounded, then $\Phi$ satisfies
condition (C); $b$ is the smallest cardinality of a family $\Phi$ of
sequences converging to $\infty$ which satisfies condition (C).
\par
{\bf Remark 2.} The family of sequences constructed in [5],
which has $b$ sequences converging to $\infty$ and satisfies
condition (C), satisfies condition (2) and does not satisfy
condition (1).
\par
{\bf Remark 3a.} ([5]) If ${\bf R}\supseteq X\not\in M$ then the
family of sequences $\{\{ x+\log (n+1)\}_{n\in\omega} : x\in X\}$
satisfies condition (C) (note: it implies well-known inequality
$b\leq\mbox{ non}(M))$ and conditions (1)-(2).
\par
{\bf Remark 3b.} Let ${\bf Q}^+$ denote the set of positive rational
numbers. Assume that: $\{ x_n \}_{n\in\omega}$ converges to $\infty$,
$\psi :\omega\longrightarrow {\bf Q}^+$ is a bijection and
$F\subseteq\omega^\omega$ satisfies $ \neg \ \exists \
g\in\omega^\omega \ \forall f\in F \ \forall^\infty k \ g(k)\neq f(k)$.
We claim that the family
$\{\{ x_n+(\psi\circ f)(n)\}_{n\in\omega}:f\in F\}$
of sequences converging to $\infty$ satisfies conditions (C) and (1).
\par
\vskip 0.3truecm
{\bf Observation.} Assume that: $\{ a_n\}$ converges to $\infty$,
$\{a_n\}$ satisfies (2) and the interval
$(k, k +\frac{1}{j}) (j\in\omega\setminus\{ 0\}, k\in\omega)$
contains more than one term of the sequence $\{ a_n\}$.
We claim that
$$
\bigcup_{n\in\omega}(a_n -\frac{1}{j}, a_n +\frac{1}{j}) \supseteq
(k,\infty)
$$
i.e. $f_{\{ a_n\} } (k)\geq j$. \\[.1cm]
\par
{\it Proof}. We choose the smallest $m\in\omega$ such that
$a_m\in (k,k +\frac{1}{j})$. Then $a_{m+1}\in (k, k +\frac{1}{j})$
and obviously $a_m -\frac{1}{j} < k$. From condition (2) we conclude
that for each $n\in\omega$, $n\geq m$:
$$
0 < a_{n+1} - a_n\leq a_n - a_{n-1}\leq\ldots\leq a_{m+2} - a_{m+1}
\leq a_{m+1} - a_m < \frac{1}{j}
$$
It implies that for each $n\in\omega$, $n > m$ we have:
$$
(a_n -\frac{1}{j}, a_n +\frac{1}{j})\supseteq (a_{n-1}, a_{n+1})
$$
Hence:
$$
\bigcup_{\begin{array}{c} n\in\omega \\ [-.2cm] n>m\end{array}}
(a_n-\frac{1}{j}, a_n +\frac{1}{j})
\ \ \supseteq\bigcup_{\begin{array}{c}
n\in\omega \\ [-.2cm] n>m\end{array}} (a_{n-1}, a_{n+1})=(a_m, \infty)
$$
Therefore:
$$
\bigcup_{n\in\omega} (a_n-\frac{1}{j}, a_n+\frac{1}{j})\supseteq
(a_m-\frac{1}{j}, a_m+\frac{1}{j})\cup (a_m, \infty)=
(a_m-\frac{1}{j}, \infty)\supseteq (k, \infty)
$$
\def\thelem{}
\begin{lem}
Assume that the family $\Phi$ of sequences converging to $\infty$
satisfies conditions (1)-(2) and $f$ is a bound of the family
$\{ f_{\{ a_n\} }:\{ a_n\}\in\Phi\}$. We claim that for each
$\{ a_n\}\in\Phi$ there exists an $m\in\omega$ such that
$\{ a_n\} $ has no more than one term in the interval
$(k, k + \frac{1}{f(k)+1})$ for each $k\in\omega$, $k\geq m$.
\end{lem}
\par
{\it Proof}. On the contrary suppose that there exist $k_0 < k_1 < k_2
<\ldots\in\omega$ such that each interval
$(k_i, k_i +\frac{1}{f(k_i)+1})$
contains more than one term of the sequence $\{ a_n\} $.
By the Observation it implies that
$$
\begin{array}{c}
f_{\{ a_n\} }(k_0)\geq f(k_0)+1 \\
f_{\{ a_n\} }(k_1)\geq f(k_1)+1 \\
f_{\{ a_n\} }(k_2)\geq f(k_2)+1 \\
..................................
\end{array}
$$
Hence $f_{\{ a_n\} }\leq_* f$ does not hold which contradicts our
assumption about $f$. This completes the proof of the Lemma.
\vspace{0.1cm}
\par
We begin the main part of the proof. Let
$f:\omega\longrightarrow\omega$ be a bound of the family
$\{ f_{\{ a_n \} }: \{ a_n \}\in\Phi\}$.
\par
For each $k\in\omega$ there exists a sequence $\{ (c(k,i),
d(k,i))\}_{i\in\omega}$ of non--empty pairwise disjoint intervals
satisfying
\begin{description}
\item{($\ast$)}
\centerline{$\bigcup_{i\in\omega} (c(k,i), d(k,i))
\subseteq (k, k+\frac{1}{f(k)+1})$}
\end{description}
We assign to each $\{ a_n\} \in\Phi$ the function $s_{\{ a_n\}}:
\omega\longrightarrow\omega$ according to the following rules:
\begin{itemize}
\item[1)] if the sequence $\{ a_n\} $ has no terms in
$\bigcup_{i\in\omega}
(c(k,i), d(k,i))$ then $s_{\{ a_n\} } (k)=0$,
\item[2)] if the sequence $\{ a_n\} $ has some term in
$\bigcup_{i\in\omega}
(c(k,i), d(k,i))$ then $s_{\{ a_n\} } (k)$ is the smallest
$i\in\omega$ such that the sequence $\{ a_n\} $ has some term in
the interval $(c(k,i), d(k,i))$.
\end{itemize}
Suppose, contrary to our claim, that the family $\Phi$ satisfies
(1)-(2) and card $\Phi < \mbox{non}(M)$. It implies that the
cardinality of the family
$\{ s_{\{ a_n\} }:\{ a_n\}\in\Phi\}\subseteq\omega^\omega$ is also
less than
$$
\begin{array}{c}
\mbox{non}(M) = \\ [.2cm]
\min \{\mbox{card } F: F\subseteq\omega^\omega
\mbox{ and } \neg \ \exists \ g\in\omega^\omega \ \forall f\in F \
\forall^\infty k \ g(k)\neq f(k)\}
\end{array}
$$
Therefore, there exists a function $g:\omega\longrightarrow\omega$
such that for each sequence
$\{ a_n\} \in\Phi \ \forall^\infty k \ g(k)\neq s_{\{ a_n\} }(k)$.
Let
$$
U:=\bigcup_{k\in\omega}(c(k,g(k)), d(k,g(k)))
$$
The set $U$ is open and unbounded above. We prove that each sequence
$\{ a_n\} \in\Phi$
has only a finite number of terms belonging to $U$. We fix
$\{ a_n\} \in\Phi$. Let $A:=\{ k\in\omega : g(k)=s_{\{ a_n\} }(k)\}$.
The set $A$ is finite. Let $B$ denote the set of all $k\in\omega$
such that the sequence $\{ a_n\} $ has more than one term in
$\bigcup_{i\in\omega}(c(k,i), d(k,i))$. The set $B$ is finite
according to the Lemma and ($\ast$). We have:
$$
U =\bigcup_{k\in A\cup B}(c(k,g(k)),
d(k,g(k)))\ \ \cup\bigcup_{k\in\omega\setminus (A\cup B)} (c(k,g(k)),
d(k,g(k)))
$$
Since $A\cup B$ is finite (as the sum of two finite sets) the first
sum over $k\in A\cup B$ is bounded above. Therefore the sequence
$\{ a_n\} $ converging to $\infty$ has only a finite number of terms
in the first sum. Concerning the second sum over
$k\in\omega\setminus (A\cup B): \mbox{ if } k\in\omega(A\cup B)$
then the sequence $\{ a_n\} $ has at most one term in
$\bigcup_{i\in\omega}(c(k,i), d(k,i))$ and this term belongs to the
interval
$(c(k, s_{\{ a_n\} }(k)), d(k, s_{\{ a_n\} }(k)))$. Since
$g(k)\neq s_{\{ a_n\} }(k)$, the interval $(c(k, g(k), d(k, g(k)))$
has no terms of the sequence $\{ a_n\} $ as it is disjoint
from $(c(k, s_{\{ a_n\} }(k)), d(k, s_{\{ a_n\} }(k)))$.
Thus $\{ a_n\} $ has no terms in
$\bigcup_{k\in\omega\setminus (A\cup B)} (c(k,g(k)),d(k,g(k)))$.
\par
We have proved that any $\{ a_n\} \in\Phi$ has only a finite number
of terms in $U$.
Therefore $\Phi$ does not satisfy condition (C). It contradicts our
assumption. The proof is complete.
\vskip 0.3truecm
\par
The following Remarks 4-5 enable simpler characterization of non$(M)$.
\vskip 0.3truecm
\par
{\bf Remark 4.} Our theorem remains valid if we replace conditions
(1)-(2) by the following condition (3) (note: condition (3) implies
condition (1)):
\begin{equation}
\mbox{ if } \{ a_n\}\in\Phi \mbox{ then }
\exists\ r>0\ \forall n\in\omega \ a_{n+1}-a_{n} > r
\end{equation}
To see this, it is sufficient to observe that condition (3) implies
the thesis of the Lemma.
\par
{\bf Remark 5.} If $( 0, \infty ) \supseteq X\not\in M$ then
the family $\Phi=\{\{(n+1)x\}_{n\in\omega} : x\in X\}$
satisfies the condition (C) ([5]), obviously $\Phi$ satisfies
condition (3).
\vskip 0.3truecm
\par
Assume that $m\in\omega\setminus\{ 0\}$.
For families $\Phi$ of sequences converging to $\infty$ we define the
following condition $(C_m)$:
\begin{itemize}
\item[($C_m$)]
{for every real sequence $x_0<y_0<x_1<y_1<x_2<y_2<...$ converging to
$\infty$ there exists a sequence $\{ a_n \} \in \Phi$
with the property that for infinitely many $k \in \omega$
the interval $(x_k, y_k)$ contains at least $m$ terms of the sequence
$\{ a_n \}$.}
\end{itemize}
\par
Obviously:
\\
$(C) \Leftrightarrow (C_1) \Leftarrow (C_2)
\Leftarrow (C_3) \Leftarrow ... \
\{f_{\{ a_n \} }: \{ a_n \}\in\Phi\}$ is unbounded
\vskip 0.3truecm
\par
Replacing $(k, k + \frac{1}{f(k)+1})$ in ($\ast$) by $(x_k, y_k)$ we
obtain the proof of the following Remark 6.
\par
{\bf Remark 6.} Our theorem remains valid if we replace conditions
(1)-(2) by the following condition (4) (note: condition (4) implies
condition (1)):
\begin{equation}
\Phi \mbox{ does not satisfy condition } (C_2)
\end{equation}
\par
{\bf Conclusion.} non($M$) is the smallest cardinality of a family
$\Phi$ of sequences converging to $\infty$ which satisfies condition
$(C_1)$ and does not satisfy condition $(C_2)$.
\vskip 1.9truecm
\centerline{\bf References}
\begin{enumerate}
\item T.~Bartoszy\'{n}ski, {\it Combinatorial aspects of
measure and category}, Fund. Math. 127 (1987), pp.225-239.
\item T.~Bartoszy\'{n}ski and H.~Judah, {\it Set theory: on the
structure of the real line}, A.~K.~Peters Ltd., Wellesley MA 1995.
\item R.~Frankiewicz and P.~Zbierski, {\it Hausdorff gaps and
limits}, North-Holland, Amsterdam 1994.
\item A.~W.~Miller, {\it A characterization of the least
cardinal for which the Baire category theorem fails}, Proc. Amer.
Math. Soc. 86 (1982), pp.498-502.
\item A.~Tyszka, {\it On a combinatorial property of families
of sequences converging to $\infty$}, J. Nat. Geom. 11 (1997),
No.1, pp.51-56.
\end{enumerate}
\begin{flushleft}
{\it Technical Faculty\\
Hugo Ko{\l}{\l}\c{a}taj University\\
Balicka 104, PL-30-149 Krak\'ow, Poland\\
rttyszka@cyf-kr.edu.pl\\
http://www.cyf-kr.edu.pl/\symbol{126}rttyszka}
\end{flushleft}
\end{document}